\pdfoutput=1
\RequirePackage{ifpdf}
\ifpdf 
\documentclass[pdftex]{sigma}
\else
\documentclass{sigma}
\fi

\input xy
\xyoption {all} \CompileMatrices

\DeclareMathOperator{\Aut}{Aut}

\newtheorem{Theorem}{Theorem}[section]
\newtheorem{Corollary}[Theorem]{Corollary}
\newtheorem{Proposition}[Theorem]{Proposition}
{\theoremstyle{definition}
\newtheorem{Example}[Theorem]{Example}
\newtheorem{Remark}[Theorem]{Remark}
}

\begin{document}

\allowdisplaybreaks

\renewcommand{\thefootnote}{$\star$}

\renewcommand{\PaperNumber}{049}

\FirstPageHeading

\ShortArticleName{The Classif\/ication of All Crossed Products $H_4 \#  k[C_{n}]$}

\ArticleName{The Classif\/ication of All Crossed Products $\boldsymbol{H_4 \#  k[C_{n}]}$\footnote{This paper is a~contribution to the Special Issue on
Noncommutative Geometry and Quantum Groups in honor of Marc A.~Rief\/fel.
The full collection is available at
\href{http://www.emis.de/journals/SIGMA/Rieffel.html}{http://www.emis.de/journals/SIGMA/Rieffel.html}}}

\Author{Ana-Loredana AGORE~$^{\dag\ddag}$, Costel-Gabriel BONTEA~$^{\S\dag}$ and Gigel MILITARU~$^\S$}

\AuthorNameForHeading{A.L.~Agore, C.G.~Bontea and G.~Militaru}

\Address{$^{\dag}$~Faculty of Engineering, Vrije Universiteit Brussel, Pleinlaan 2, B-1050 Brussels, Belgium}
\EmailD{\href{mailto:ana.agore@vub.ac.be}{ana.agore@vub.ac.be}, \href{mailto:ana.agore@gmail.com}{ana.agore@gmail.com}}
\URLaddressD{\url{http://homepages.vub.ac.be/~aagore/}}

\Address{$^{\ddag}$~Department of Applied Mathematics, Bucharest University of Economic Studies,\\
\hphantom{$^{\ddag}$}~Piata Romana 6, RO-010374 Bucharest 1, Romania}

\Address{$^{\S}$~Faculty of Mathematics and Computer Science,
University of Bucharest,\\
\hphantom{$^{\S}$}~Str.~Academiei 14, RO-010014 Bucharest~1, Romania}
\EmailD{\href{mailto:costel.bontea@gmail.com}{costel.bontea@gmail.com}, \href{mailto:gigel.militaru@fmi.unibuc.ro}{gigel.militaru@fmi.unibuc.ro},
\href{mailto:gigel.militaru@gmail.com}{gigel.militaru@gmail.com}}
\URLaddressD{\url{http://fmi.unibuc.ro/ro/departamente/matematica/militaru_gigel/}}

\ArticleDates{Received November 18, 2013, in f\/inal form April 18, 2014; Published online April 23, 2014}

\Abstract{Using the computational approach introduced in
[Agore~A.L., Bontea~C.G., Militaru~G., \textit{J.~Algebra Appl.} \textbf{12} (2013), 1250227,
  24~pages] we classify all coalgebra split extensions of $H_4$ by
$k[C_n]$, where $C_n$ is the cyclic group of order $n$ and $H_4$
is Sweedler's $4$-dimensional Hopf algebra. Equivalently, we
classify all crossed products of Hopf algebras $H_4 \#
k[C_{n}]$ by explicitly computing two classifying objects: the
cohomological `group' ${\mathcal H}^{2} ( k[C_{n}], H_4)$ and
$\textsc{C}\textsc{r}\textsc{p} ( k[C_{n}], H_4):=$ the set of
types of isomorphisms of all crossed pro\-ducts~$H_4 \#
k[C_{n}]$. More precisely, all crossed products $H_4 \# k[C_n]$
are described by generators and relations and classif\/ied: they are
$4n$-dimensional quantum groups $H_{4n,  \lambda,  t}$,
parame\-teri\-zed by the set of all pairs $(\lambda,  t)$ consisting
of an arbitrary unitary map $t : C_n \to C_2$ and an $n$-th root~$\lambda$ of~$\pm 1$. As an application, the group of Hopf algebra
automorphisms of $H_{4n,  \lambda,  t}$ is explicitly
described.}

\Keywords{crossed product of Hopf algebras; split extension of Hopf algebras}

\Classification{16T10; 16T05; 16S40}

\rightline{\it Dedicated to Marc Rieffel on the occasion of his 75th birthday}

\renewcommand{\thefootnote}{\arabic{footnote}}
\setcounter{footnote}{0}

\section{Introduction}

The second cohomology group ${\rm H}^2 (H, A)$ classif\/ies all
extensions of an Abelian group~$A$ by a~group~$H$, i.e.\ all groups~$G$ that f\/it into an exact sequence $1 \to A \to G \to H \to 1$.
More precisely, each element $f \in {\rm H}^2 (H, A)$ is assigned
with an extension~$G_{f}$ of~$A$ by~$H$, namely the crossed
product $G_f := A \#^f H$ of~$A$ and~$H$ and the ef\/fective
classif\/ication of all the extensions of~$A$ by~$H$ is obtain after
computing the group ${\rm H}^2 (H, A)$. In group theory there is a
well developed cohomological machinery \cite{adem} which allows,
at least for some classes of groups~$A$ and~$H$, to compute this
cohomology group. Transferring this problem to Hopf algebras by
considering group algebras over a f\/ield~$k$, we obtain that any
extension of~$A$ by~$H$ gives a coalgebra split extension of~$k[A]$ by~$k[H]$ in the sense of~\cite[Def\/inition~1.2]{ABM2}. In
fact, there is more: a~Hopf algebra~$E$ is a coalgebra split
extension of~$k[A]$ by~$k[H]$ if and only if $E \cong k[G]$, for a
group~$G$ which is an extension of $A$ by $H$ \cite[Example~1.4]{ABM2}. This can be restated in a cohomological manner as
follows: ${\mathcal H}^{2} (k[H], k[A]) \cong {\rm H}^{2} (H, A)$,
where ${\mathcal H}^{2} (k[H], k[A])$ denotes the second
cohomological group for Hopf algebras introduced by Sweedler~\cite{Sw68}. Thus, the classif\/ication of all coalgebra split
extensions of a Hopf algebra~$A$ by a Hopf algebra $H$ covers the
extension problem from group theory. Now, if we replace the group
algebra~$k[A]$ by another arbitrary Hopf algebra, say a
noncommutative and noncocommutative Hopf algebra such as
Sweedler's $4$-dimensional Hopf algebra~$H_4$, things change
radically as none of the classical cohomological techniques can be
applied in this context. In other words, the cohomological type
object ${\mathcal H}^{2} (k[H], H_4)$ needs to be computed using a
direct approach.

In this paper we shall classify all coalgebra split extensions of
$H_4$ by $k[C_n]$, where $C_n$ is the cyclic group of order $n$,
i.e.\ all Hopf algebras $E$ that f\/it into a sequence $H_4
\hookrightarrow E \stackrel{\pi} \to k[C_{n}] $ such that $\pi: E
\to k[C_{n}]$ splits as a coalgebra map and $H_4 \simeq E^{{\rm
co}(k[C_{n}])}$. Equivalently, we classify all crossed products of
Hopf algebras $H_4 \#  k[C_{n}]$. This kind of crossed products
are a special case of those arising in Hopf--Galois extensions
theory. For a generalization of the notion of a~Hopf--Galois
extension see~\cite{BH}. The fact that there is no ef\/f\/icient
cohomology theory for arbitrary Hopf algebras led us to consider a
purely computational approach which relies heavily on the methods
introduced in~\cite{ABM2}. First of all we compute the set of all
crossed systems $(H_{4}, k[C_{n}], \triangleright, f)$ between~$H_4$ and~$k[C_n]$. This is the f\/irst computational part of our
strategy, quite laborious considering the large number of axioms
that need to be fulf\/illed by the pairs $(\triangleright, f)$ in
order to provide a crossed system of Hopf algebras $(H_{4},
k[C_{n}], \triangleright, f)$.  Theorem~\ref{sistemecros} gives the
following description: the set of all crossed systems of Hopf
algebras $(H_{4}, k[C_{n}], \triangleright, f)$ is parameterized
by the set ${\mathcal C} {\mathcal S} (n, k) \subseteq {\mathcal
U} (C_n, C_2) \times k^*$, consisting of all pairs $(t, \lambda)$,
where $t : C_{n} \to C_{2}$ is an arbitrary unitary map and
$\lambda \in k^*$ is an $n$-th root of $\pm 1$. Thus there are at
most $n 2^{n-1}$ crossed products of the form $H_4 \# k[C_n] :=
H_{4n, \lambda,  t}$, for some $(t, \lambda) \in {\mathcal C}
{\mathcal S} (n, k) $ and they are described by generators and
relations in Corollary~\ref{anavers}~-- this is the second step of our
approach. Finally, the last computational step uses \cite[Theorem~2.1]{ABM2} as a tool: we shall classify all the above crossed
products by explicitly computing the classifying objects
${\mathcal H}^{2} ( k[C_{n}], H_4)$ and
$\textsc{C}\textsc{r}\textsc{p} ( k[C_{n}], H_4) :=$ the set of
types of isomorphisms of Hopf algebras of all crossed products
$H_4 \#  k[C_n]$. The classif\/ication results are proven in Theorems~\ref{clasificareH_4C_n} and~\ref{aldoileaseint}. As an
application, Corollary~\ref{autom} provides the parametrization of
$\Aut_{\rm Hopf} (H_{4n, \lambda,  t})$, the group of Hopf
algebra automorphisms of~$H_{4n, \lambda,  t}$.

We point out that a coalgebra split extension is a special case of
a more general type of Hopf algebra extension def\/ined in
\cite[Def\/inition~1.2.0]{AD}. As explained in~\cite[Section~5.2]{AndCan}, it is a very dif\/f\/icult task to classify this general
type of extensions and, to the best of our knowledge, the only
example of such a classif\/ication is~\cite[Lemma~2.8]{GV}. There
are also several known classif\/ication results for Abelian
extensions of Hopf algebras associated to matched pairs of groups
(see \cite{krop, mas} and the references therein). Our
computational method, introduced in \cite{ABM2} and used in the
present paper, is the most direct and natural way of approaching
the classif\/ication of all coalgebra split extensions for two given
Hopf algebras. This might not be the only way to approach the
problem: for certain pairs of Hopf algebras~$A$ and~$H$, the class
of all coalgebra split extensions of~$A$ by~$H$ can be classif\/ied
by using the lifting method~\cite{AS1}.

\section{Preliminaries} \label{preli}
Let $G$ and $H$ be two groups with $H$ Abelian. In what follows
$\text{ord} (g)$ denotes the order of the element $g \in G$, $|G|$
is the order of $G$ while $\Aut (G)$ stands for the group of
automorphisms of $G$. $\mathcal{U} (G, H)$ will be the set of all
unitary maps $u: G \to H$, i.e.~$u(1) = 1$. \emph{A normalized
$2$-cocycle} \cite[Section~7]{mont} is a map $f : G \times G \to
H$ such that
\begin{gather}\label{standcoc}
f(1,  g) = f(g, 1) = 1, \qquad f(g,  h) f(gh,  l) = f(h,
 l) f(g,  hl)
\end{gather}
for all $g, h, l \in G$. The map $f : G \times G \to H$, $f
(g, h) :=1$, for all $g,h\in G$, is called the trivial cocycle.
The set of all normalized $2$-cocycles is denoted by $Z^2 (G, H)$.
Two cocycles~$f$ and~$f'$ are called \emph{cohomologous}, and we
denote this by $f \approx f'$, if there exists a unitary map $u: G
\to H$ such that
\begin{gather*}
f (g, g') = u(g)  u(g')  f' (g, g')  u (gg')^{-1}
\end{gather*}
for all $g,g'\in G$. A normalized $2$-cocycle $f$ is called a
\emph{coboundary} if $f$ is cohomologous with the trivial
cocylcle, i.e.\ if there exists a unitary map $u: G \to H$ such
that $f (g, g') = f_u (g, g') := u(g) u(g') u (gg')^{-1}$, for all
$g$, $g'\in G$. ``$\approx$'' is an equivalence relation on $Z^2 (G,
H)$ and we denote by ${\mathcal H}^2 (G, H) := Z^2 (G, H)/{\approx}$ the corresponding factor set. ${\mathcal H}^2 (G, H)$ is called
the second cohomology group of $G$ with coef\/f\/icients in the
Abelian group~$H$~\cite{adem}. The following result is folklore:
it describes all the normalized $2$-cocycles $f : C_n \times C_n
\to C_2$.

\begin{Proposition}\label{prop1}
Let $n$ be a positive integer. There exists a bijection ${\mathcal
U} (C_n, C_2) \simeq Z^2 (C_n, C_2)$ such that the $2$-cocycle $f
= f_t : C_n\times C_n \to C_2$ associated to the unitary map $t
\in \mathcal{U} (C_n, C_2) $ is given by
\begin{gather*}
f_t \big(c^{i}, c^{j}\big) := \prod_{p=0}^{i-1} t (c^p) \prod_{q=0}^{j-1}
t (c^q) \prod_{r=0}^{i+j-1} t (c^r)
\end{gather*}
for all $i, j = 1, \dots, n$.
\end{Proposition}

\begin{proof}
The fact that $f_t$ is a normalized $2$-cocycle follows from a
straightforward computation and moreover, $f_t$ can be written
equivalently as
\begin{gather}\label{cocciliceb}
f_t\big(c^{i}, c^{j}\big) = t\big(c^0\big) t\big(c^1\big) \cdots t \big(c^{i-1}\big)   t\big(c^j\big)
t\big(c^{j+1}\big) \cdots t \big(c^{j+i-1}\big)
\end{gather}
for all $i, j = 1, \dots, n$. The inverse of the map $t \mapsto
f_t$ is constructed as follows: let $f: C_n \times C_n \to C_2$ be
a normalized $2$-cocycle and def\/ine $t = t_f: C_n \to C_2$, $t
(c^j) := f (c,  c^j)$, for all $j = 1, \dots, n$. Using the
cocycle condition~\eqref{standcoc} and induction on $i$ we easily
obtain that~\eqref{cocciliceb} holds and the correspondence $(f
\mapsto t_f, t\mapsto f_t)$ is bijective.
\end{proof}

\subsection*{Crossed products of Hopf algebras}

We shall review the construction of the crossed product of two
Hopf algebras introduced in \cite[Lemma~1.2.10]{AndN} as a special
case of the cocycle bicrossproduct \cite[Theorem~2.20]{AD}. It can
be also obtained as a special case of the unif\/ied product of
\cite[Theorem~2.4, Examples~2.5(2)]{am1}. From now on~$k$ will be
an arbitrary f\/ield and we shall use $\otimes$ instead of~$\otimes_{k}$. For the comultiplication of a~Hopf algebra we use
Sweedler's $\Sigma$-notation with suppressed summation sign:
$\Delta (a) = a_{(1)}\otimes a_{(2)}$. If~$A$ and~$H$ are Hopf
algebras and $f: H \otimes H \to A$ is a $k$-linear map, we denote
$f (g, h) = f (g\otimes h)$, for all $g,h\in H$. A $k$-linear
map $\triangleright: H \otimes A \to A$ is called a \textit{weak
action}~\cite{mont} of $H$ on $A$ if for any $a,b \in A$, $h
\in H$:
\begin{gather}
 1_H \triangleright a  =  a,  \label{1b} \\
 h \triangleright 1_{A}  =  \varepsilon_{H} (h) 1_{A}, \label{2}
\\
 h \triangleright (ab)  =  (h_{(1)} \triangleright a)
(h_{(2)} \triangleright b). \label{3}
\end{gather}
Let $A$ and $H$ be two Hopf algebras, $f: H \otimes H \to A$ a
$k$-linear map and $\triangleright: H \otimes A \to A$ a weak
action of~$H$ on~$A$. We denote by $A \#_{f}^{\triangleright} H$
the $k$-vector space $A \otimes H$ endowed with the following
multiplication
\begin{gather}\label{multiplicarea}
(a \# g) \cdot (b \# h):= a (g_{(1)} \triangleright b) f \big(
g_{(2)} , h_{(1)} \big) \#  g_{(3)} h_{(2)}
\end{gather}
for all $a, b \in A$, $g, h \in H$, where $\#$ stands for
$\otimes$. The object $A \#_{f}^{\triangleright} H$ is called the
\textit{crossed product} of $A$ with $H$ if it is a~Hopf algebra
with the multiplication~\eqref{multiplicarea}, the unit $1_A \#
1_H$ and the coalgebra structure given by the tensor product of
coalgebras. \cite[Proposition~1.1]{ABM2} and \cite[Lemma~1.2.10]{AndN} proves that $A \#_{f}^{\triangleright} H$ is a~crossed product of $A$ with $H$ if and only if~$f$ and~$\triangleright$ are morphisms of coalgebras satisfying the
following compatibilities
\begin{gather}
f (h, 1_{H})  =  f (1_{H}, h) = \varepsilon_{H}(h)
1_{A}, \label{1a} \\
\big[g_{(1)} \triangleright (h_{(1)} \triangleright a) \big] f
(g_{(2)},
h_{(2)})  =  f (g_{(1)}, h_{(1)}) \big((g_{(2)} h_{(2)}) \triangleright a \big), \label{4} \\
\big( g_{(1)} \triangleright f (h_{(1)}, l_{(1)}) \big) f \big(
g_{(2)}, h_{(2)} l_{(2)} \big)  =
f (g_{(1)}, h_{(1)}) f (g_{(2)}h_{(2)}, l), \label{5}\\
g_{(1)} \otimes ( g_{(2)} \triangleright a)  =
g_{(2)} \otimes (g_{(1)} \triangleright a), \label{6}\\
g_{(1)} h_{(1)} \otimes f (g_{(2)}, h_{(2)})  =  g_{(2)} h_{(2)}
\otimes f(g_{(1)}, h_{(1)}) \label{7}
\end{gather}
for all $a \in A$, $g, h, l \in H$. Whenever
\eqref{1b}--\eqref{3} and~\eqref{1a}--\eqref{7} hold, $(A, H,
\triangleright, f)$ is called a \textit{crossed system of Hopf
algebras}. The antipode of $A \#_{f}^{\triangleright}   H$ is
given for any $a \in A$ and $h \in H$ by
\begin{gather}\label{antipod}
S(a \# h) := \big( S_{A} \big[ f \big( S_{H} (h_{(2)}), h_{(3)}
\big) \big] \# S_{H} (h_{(1)}) \big) \cdot \big( S_{A} (a) \#
1_{H} \big).
\end{gather}
A \emph{coalgebra split extension of $A$ by $H$} \cite{ABM2} is a
triple $(E, i, \pi)$ consisting of a Hopf algebra $E$ and two Hopf
algebra homomorphisms $i: A \to E$ and $\pi : E \to H$ such that
$i$ is injective, $\pi$ has a section as a coalgebra map and $i(A)
\cong E ^{{\rm co}(H)} := \{x \in E \, | \, x_{(1)} \otimes \pi
(x_{(2)}) = x \otimes 1 \} $. If $E$ is f\/inite-dimensional, then
the sequence $\xymatrix{ A \ar[r]^{i} & {E} \ar[r]^{\pi} & H}$ is
an exact sequence of Hopf algebras~\cite{AD}. Two coalgebra split
extensions $(E, i, \pi)$, $(E', i', \pi')$ of $A$ by $H$ are
called \emph{equivalent} if there exists an isomorphism of Hopf
algebras $\psi : E \to E'$ that stabilizes $A$ and co-stabilizes~$H$, i.e.\ the following diagram
\begin{gather*}
\xymatrix {& A \ar[r]^{i} \ar[d]_{{\rm Id}_{A}} & {E}
\ar[r]^{\pi}\ar[d]^{\psi} & H\ar[d]^{{\rm Id}_{H}}\\
& A\ar[r]^{i'} & {E'}\ar[r]^{\pi'} & H}
\end{gather*}
is commutative. Any crossed product $A \#_{f}^{\triangleright}
H$ is a coalgebra split extension of $A$ by $H$ via $i_A : A \to A
\#^{\triangleright}_{f} H$, $i_A (a) = a  \# 1_H$, for all $a \in
A$, and $\pi_H : A \#^{\triangleright}_{f} H \to H$, $\pi_H (a \#
h) = \varepsilon_A (a) h$, for all $a \in A$ and $h \in H$.
Conversely, any coalgebra split extension of~$A$ by a
cocommutative Hopf algebra $H$ is equivalent to a crossed product
extension $(A \#_{f}^{\triangleright} H, i_{A}, \pi_{H})$ of~$A$
by $H$ \cite[Proposition~1.3]{ABM2}. More precisely, if $\varphi:
H \to E$ is a unit preserving coalgebra map that splits~$\pi$ then
the action $\triangleright = \triangleright_{\varphi}$ and the
cocycle $f = f_{\varphi}$ implemented by $\varphi$ are given by:
$h \triangleright a := \varphi (h_{(1)})   a   S \varphi
(h_{(2)})$ and $f (g,   h) := \varphi (g_{(1)})   \varphi
(h_{(1)})   S \varphi (g_{(2)} h_{(2)})$, for all $g,h\in H$
and $a\in A$. Thus, the classif\/ication of all coalgebra split
extensions of $A$ by $H$ reduces to the classif\/ication of all
crossed products $A \#_{f}^{\triangleright} H$. The classifying
object for all coalgebra split extensions of $A$ by $H$, denoted
by $\mathcal{H}^{2} (H, A)$, was introduced in \cite[Remark~2.4]{ABM2}. For the reader's convenience we recall brief\/ly its
construction: let ${\mathcal C} {\mathcal S} (A, H)$ be the set of
all pairs $(\triangleright, f)$ such that $(A, H, \triangleright,
f)$ is a crossed system of Hopf algebras. Two pairs
$(\triangleright, f)$ and $(\triangleright', f') \in {\mathcal
C}{\mathcal S} (A, H)$ are called \emph{cohomologous} and we
denote this by $(\triangleright, f) \approx (\triangleright', f')$
if there exists $r: H \rightarrow A$ an unitary cocentral
map\footnote{I.e.\ $r$ is a unit preserving coalgebra map such that
$r(h_{(1)}) \otimes h_{(2)} = r(h_{(2)}) \otimes h_{(1)}$, for all
$h\in H$.} such that for any $a\in A$ and $h,g \in H$ we have
\begin{gather*}
h \triangleright' a  =  r (h_{(1)})   (h_{(2)}
\triangleright a)   (S_A \circ r )(h_{(3)}), \\
f'(h,   g)  =  r(h_{(1)})   \big(h_{(2)} \triangleright
r(g_{(1)}) \big)   f(h_{(3)},   g_{(2)})    (S_A \circ r)
(h_{(4)} g_{(3)}).
\end{gather*}

Then \cite[Theorem~3.4]{am1} proves that $(\triangleright, f)
\approx (\triangleright', f')$ if and only if there exists a Hopf
algebra isomorphism $A \#_{f}^{\triangleright}  H \cong A
\#_{f'}^{\triangleright'}  H$ that stabilizes~$A$ and
co-stabilizes~$H$. Thus, $\approx$ is an equivalence relation on
the set ${\mathcal C} {\mathcal S} (A, H)$ and $ {\mathcal H}^{2}
(H, A)$ is the pointed quotient set def\/ined by ${\mathcal H}^{2}
(H, A) := {\mathcal C} {\mathcal S} (A, H)/ {\approx}$. Now, if $H$
is cocommutative and $A$ commutative then ${\mathcal H}^{2} (H,
A)$ coincides with the second cohomological group as constructed
by Sweedler~\cite{Sw68}; however, for arbitrary Hopf algebras~$A$
and~$H$, we could not f\/ind a complex for which ${\mathcal H}^{2}
(H, A)$ is the associated cohomological group. We denote by
$\textsc{C}\textsc{r}\textsc{p} (H, A)$ the set of types of Hopf
algebra isomorphisms of all crossed products $A
\#_{f}^{\triangleright}   H$ associated to all crossed systems
$(A, H, \triangleright, f)$. Two equivalent extensions are
isomorphic and hence there exists a canonical surjection
${\mathcal H}^{2} (H, A) \twoheadrightarrow
\textsc{C}\textsc{r}\textsc{p} (H, A)$. As in the case of groups,
it turns out that the two classifying objects ${\mathcal H}^{2}
(H, A)$ and $\textsc{C}\textsc{r}\textsc{p} (H, A)$ are dif\/ferent:
\cite[Proposition~4.2]{ABM2} proves that for the pair of Hopf
algebras $(H, A) := (H_4, k[Y])$ we have that
$\textsc{C}\textsc{r}\textsc{p} (H, A) \cong C_2$ while ${\mathcal
H}^{2} (H, A) \cong k$, where~$k$ is the base f\/ield.

\section[The classif\/ication of coalgebra split extensions of $H_{4}$ by $\protect{k[C_{n}]}$]{The classif\/ication of coalgebra split extensions of $\boldsymbol{H_{4}}$ by $\boldsymbol{k[C_{n}]}$}
\label{exclasros}

Let $k$ be a f\/ield of characteristic $\ne 2$, $k[C_{n}]$ the group
Hopf algebra of the cyclic group $C_{n}$ of order $n$ generated by
$c$ and let $H_{4}$ be Sweedler's 4-dimensional Hopf algebra
generated by the group-like element $g$ and the $(1, g)$-primitive
element $x$ (that is $\Delta (x) = x \otimes 1 + g \otimes x$)
subject to the following relations
\begin{gather*}
g^{2} = 1, \qquad x^{2} = 0, \qquad x g = -g x.
\end{gather*}
We classify all coalgebra split extensions of $H_{4}$ by
$k[C_{n}]$ following the strategy presented in the introduction:
f\/irst, we describe explicitly all crossed systems $(H_{4},
k[C_{n}], \triangleright, f)$, then we describe by generators and
relations the associated crossed products and, f\/inally, we
classify such extensions by computing the classifying objects
${\mathcal H}^{2} ( k[C_{n}], H_4)$ and
$\textsc{C}\textsc{r}\textsc{p} ( k[C_{n}], H_4)$. To any unitary
map $t : C_{n} \to C_{2}$ we associate the function $\sigma_{t} :
\mathbb{N} \to \mathbb{N}$ def\/ined by $\sigma_{t}(0) =
\sigma_{t}(1) := 0$ and for any positive integer  $j \geq 2$
\begin{gather*}
 \sigma_{t}(j) :=
\begin{cases}
1 + \textnormal{ord} \, t(c) + \cdots + \textnormal{ord} \,
t\big(c^{j-1}\big) & \textnormal{if}  \quad 2 \mid j,\\
\textnormal{ord} \, t(c) + \cdots + \textnormal{ord} \, t\big(c^{j-1}\big)
& \textnormal{if} \quad  2 \nmid j.
\end{cases}
\end{gather*}

The following gives the parametrization of all crossed systems
$(H_{4}, k[C_{n}], \triangleright, f)$.

\begin{Theorem} \label{sistemecros}
Let $k$ be a field of characteristic ${\neq}2$ and $n$ a positive
integer. Then there exists a~bijective correspondence between the
set of all crossed systems of Hopf algebras $(H_4, k[C_n],
\triangleright, f)$ and the set $\mathcal{CS} (n, k) \subseteq
\mathcal{U} (C_n, C_2) \times k^*$ consisting of all pairs $(t,
\lambda)$, where $t : C_{n} \to C_{2}$ is a~unitary map and
$\lambda \in k^*$ is such that $\lambda^{n} =
(-1)^{\sigma_{t}(n)}$.

Under this bijection the crossed system $(H_4, k[C_n],
\triangleright, f)$ corresponding to $(t, \lambda)$ is given by
\begin{gather}
c^{j} \triangleright 1 = 1, \qquad c^{j} \triangleright g = g,
\qquad c^{j} \triangleright x = (-1)^{\sigma_{t}(j)} \lambda^{j}
  x, \qquad c^{j} \triangleright (gx) = (-1)^{\sigma_{t}(j)}
\lambda^{j}   gx,\nonumber 
\\
f\big(c^{i}, c^{j}\big) = \prod_{p = 0}^{i-1} t(c^{p}) \prod_{q = 0}^{j-1}
t(c^{q}) \prod_{r = 0}^{i + j - 1} t(c^{r}) \label{cociclas11}
\end{gather}
for all $i, j = 1, \dots, n$.
\end{Theorem}

\begin{proof}
We have to describe all pairs $(\triangleright, f)$ of coalgebra
maps $\triangleright : k[C_n] \otimes H_4 \to H_4$ and $f : k[C_n]
\otimes k[C_n] \to H_4$ that satisfy the compatibility conditions~\eqref{1b}--\eqref{3} and~\eqref{1a}--\eqref{7}. For such a
pair, the compatibilities~\eqref{6}--\eqref{7} hold since
$k[C_n]$ is cocommutative.

Let $(\triangleright, f)$ be a pair such that $(H_4, k[C_n],
\triangleright, f)$ is a crossed system. We will f\/irst prove that
$f (c^{i}, c^{j}) \in \{1,g\}$, for all $i,j = 1, \dots, n$,
the corresponding map $f : C_n \times C_n \to C_2 = \{ 1, g \}$ is
a~classical normalized $2$-cocycle of groups in the sense of~\eqref{standcoc} and the action $\triangleright : k[C_n] \otimes
H_4 \to H_4$ is given by
\begin{gather}
c^{j} \triangleright 1 = 1, \qquad c^{j} \triangleright g = g,
\qquad c^{j} \triangleright x = \lambda_{j}  x, \qquad c^{j}
\triangleright (gx) = \lambda_{j}   gx \label{actiunea_rhd}
\end{gather}
for all $j = 1, \dots, n$, where the scalars $\lambda_{j} \in k$
are such that \eqref{actiunea_rhd} 
 holds, which in our case
takes the equivalent form
\begin{gather}
\big( c^{i} \triangleright (c^{j} \triangleright a) \big)
f\big(c^{i},c^{j}\big) = f\big(c^{i},c^{j}\big) \big(c^{i+j} \triangleright a\big)
\label{4HC}
\end{gather}
for all $i,j = 1, \dots, n$ and $a \in H_4$. In the next step
we will use Proposition~\ref{prop1}: any classical normalized $2$-cocycle $f:
C_n \times C_n \to C_2$ is implemented by a unique unitary map $t
: C_{n} \to C_{2}$ such that~\eqref{cociclas11} holds. Finally,
using this description of~$f$, we will prove that~\eqref{4HC}
holds if and only if $\lambda_{j} = (-1)^{\sigma_{t}(j)}
\lambda_{1}^{j}$, for all $j = 1, \dots, n$, and $\lambda_{1}^{n}
= (-1)^{\sigma_{t}(n)}$. This will f\/inish the proof.

First, since $f: k[C_n] \otimes k[C_n] \to H_4$ is a coalgebra
map, $f (c^{i}, c^{j})$ is a grouplike element in $H_4$. Therefore
$f (c^{i}, c^{j}) \in \{1, g\}$, for all $i,j = 1, \dots, n$.
Hence, $f$ is uniquely determined by its restriction to~$C_{n}$,
which we will also denote by $f: C_n \times C_n \to C_2$.
Similarly, as $\triangleright : k[C_n] \otimes H_4 \to H_4$ is a
coalgebra map, we obtain that $c^{j} \triangleright g \in
\{1,g\}$, for all $j = 1, \dots, n$. We claim that $c^{j}
\triangleright g = g$, for all $j = 1, \dots, n$. Indeed, assume
that $c^{j} \triangleright g = 1$, for some $j = 1, \dots, n$.
Applying~\eqref{4HC} for $i := n-j$ and $a := g$ and taking into
account the normalizing condition~\eqref{1a}, we obtain
$f(c^{n-j}, c^{j}) = f(c^{n-j}, c^{j})  g$, which gives a
contradiction as ${\rm Im} (f) \subseteq \{1, g\}$ and $g^2 = 1$.
Therefore, $c^{j} \triangleright g = g$, for all $j = 1, \dots,
n$.

Using once again that $\triangleright$ is a coalgebra map and the
fact that $x$ is a $(1, g)$-primitive element of~$H_4$ we obtain
that $c^{j} \triangleright x$ is also an $(1, g)$-primitive
element of~$H_4$. Hence, using \cite[Lemma~4.2]{ABM1}, we obtain
that $c^{j} \triangleright x = \mu_{j} - \mu_{j} g + \lambda_{j}
x$, for some scalars $\mu_{j}$, $\lambda_{j} \in k$. From this
formula of $c^{j} \triangleright x$ and the compatibility
condition~\eqref{3} we obtain
\begin{gather*}
0 = c^ j \triangleright \big(x^2\big) = \big(c^{j} \triangleright x\big) \big(c^{j}
\triangleright x\big) = 2 \mu_{j}^{2} - 2\mu_{j}^{2} g + 2 \mu_{j}
\lambda_{j}x.
\end{gather*}
Thus $\mu_{j} = 0$ and $c^{j} \triangleright x = \lambda_{j}x$,
for all $j = 1, \dots, n$. Applying~\eqref{3} once again we
obtain: $c^{j} \triangleright (gx) = (c^{j} \triangleright g)
(c^{j} \triangleright x) = \lambda_{j} gx$. Thus, the formula~\eqref{actiunea_rhd} holds. Now, the compatibility condition~\eqref{5} is equivalent to
\begin{gather*}
\big( c^{i} \triangleright f\big(c^{j},c^{k}\big) \big)
f\big(c^{i},c^{j+k}\big) = f\big(c^{i},c^{j}\big) f\big(c^{i+j}, c^{k}\big)
\end{gather*}
for all $i, j, k = 1, \dots, n$. Since $c^{i} \triangleright f
(c^{j}, c^{k}) = f(c^{j},c^{k})$, for all $i$, $j$ and $k$, the
condition~\eqref{5} is, thus, equivalent to the fact that $f: C_n
\times C_n \to C_2$ is a usual normalized $2$-cocycle for groups
in the sense of~\eqref{standcoc}. It follows then from Proposition~\ref{prop1} that there exists a unique map $t : C_n \to C_2$ such
that $t(1) = 1$ and \eqref{cociclas11} holds
for all $i, j = 1, \dots, n$. In particular, $f(c, c^{j}) =
t(c^{j})$, for all $j = 1, \dots, n$.

With the information collected so far on $f$ and $\triangleright$
we turn now to the compatibility condition~\eqref{4HC} and see
when it is satisf\/ied: this is in fact the last compatibility that
needs to be fulf\/illed in order for $(H_4, k[C_n], \triangleright,
f)$ to be a crossed system of Hopf algebras. We notice that~\eqref{4HC} holds automatically if $a = g$, since $c^i
\triangleright g = g$, for all $i = 1, \dots, n$. On the other
hand, using~\eqref{4HC} for $a = x$ and $c^i = c^j = c$ we obtain
\begin{gather*}
\lambda_{1}^{2} x f(c,c) = f(c,c) \lambda_{2} x.
\end{gather*}
Multiplying to the left by $f(c,c) = t(c) \in \{1, g\}$ we obtain
\begin{gather*}
\lambda_{2} x = \lambda_{1}^{2} f(c,c) x f(c,c) = (-1)^{1+
\textnormal{ord} \, f(c,c)} \lambda_{1}^{2} x.
\end{gather*}
Thus, $\lambda_{2} = (-1)^{\sigma_{t}(2)} \lambda_{1}^{2}$. By
writing~\eqref{4HC} for $(c, c^2)$ and $a = x$ we obtain that
$\lambda_{1} \lambda_{2} x f(c,c^{2}) = f(c,c^{2}) \lambda_{3} x$,
from which we deduce
\begin{gather*}
\lambda_{3} x = \lambda_{1} \lambda_2 f\big(c,c^2\big) x f\big(c,c^{2}\big) =
(-1)^{\textnormal{ord} \, f(c,c) + \textnormal{ord} \, f(c,c^{2})}
\lambda_{1}^{3} x.
\end{gather*}
Thus, $\lambda_{3} = (-1)^{\sigma_{t}(3)} \lambda_{1}^{3}$. By
induction, it follows that $\lambda_{j} =
(-1)^{\sigma_{t}(j)}\lambda_{1}^{j}$, for all $j = 2, \dots, n-1$.
Applying~\eqref{4HC} for $(c, c^{n-1})$ and $a = x$, we obtain
that $\lambda_{1} \lambda_{n-1} x f(c, c^{n-1}) = f(c, c^{n-1} ) x
$, which is equivalent to
\begin{gather*}
x = \lambda_{1} \lambda_{n-1} f\big(c, c^{n-1}\big) x f\big(c, c^{n-1}\big) =
(-1)^{\sigma_{t}(n)} \lambda_{1}^{n} x.
\end{gather*}
Thus, $\lambda_{1}^{n} = (-1)^{\sigma_{t}(n)}$, as needed.

Conversely, let $\lambda \in k$ be such that $\lambda^{n} =
(-1)^{\sigma_{t}(n)}$ and def\/ine $\lambda_{j} :=
(-1)^{\sigma_{t}(j)} \lambda^{j}$. Then we can prove that~\eqref{4HC} holds. Indeed, it suf\/f\/ices to see that~\eqref{4HC}
holds for $a = x$. Observe f\/irst that, since
$(-1)^{\sigma_{t}(n+j)} = (-1)^{\sigma_{t}(n) + \sigma_{t}(j)}$,
we have $c^{n+j} \triangleright x = (-1)^{\sigma_{t}(n+j)}
\lambda^{n+j} x$, and thus, $c^{j} \triangleright x =
(-1)^{\sigma_{t}(j)} \lambda^{j} x$, for all non-negative integers
$j$. Secondly,
\begin{gather*}
x  t(c)   t\big(c^{2}\big) \cdots t\big(c^{j-1}\big)   = (-1)^{1 +
\textnormal{ord} \, t(c)} t(c)   x   t\big(c^{2}\big) \cdots t\big(c^{j-1}\big) =\cdots
\\
\hphantom{x  t(c)   t\big(c^{2}\big) \cdots t\big(c^{j-1}\big)}{}
 = (-1)^{j-1 + \textnormal{ord} \, t(c) + \cdots +
\textnormal{ord} \, t(c^{j-1})} t(c)   t\big(c^{2}\big) \cdots t\big(c^{j-1}\big)
x   \\
\hphantom{x  t(c)   t\big(c^{2}\big) \cdots t\big(c^{j-1}\big)}{}
= (-1)^{\sigma_{t}(j)} t(c)   t\big(c^{2}\big) \cdots t\big(c^{j-1}\big) x.
\end{gather*}
Using this formula, we have
\begin{gather*}
\big( c^{i} \triangleright \big(c^{j} \triangleright x\big) \big) f\big(c^{i},
c^{j}\big)   = (-1)^{\sigma_{t}(i) + \sigma_{t}(j)} \lambda^{i+j} x
f\big(c^{i},
c^{j}\big)\\
\hphantom{\big( c^{i} \triangleright \big(c^{j} \triangleright x\big) \big) f\big(c^{i},
c^{j}\big)}{}
 = (-1)^{\sigma_{t}(i) + \sigma_{t}(j)} \lambda^{i+j} x \prod_{p
= 0}^{i-1} t(c^{p}) \prod_{q = 0}^{j-1} t(c^{q}) \prod_{r = 0}^{i
+ j - 1} t(c^{r}) \\
\hphantom{\big( c^{i} \triangleright \big(c^{j} \triangleright x\big) \big) f\big(c^{i},
c^{j}\big)}{}
  = (-1)^{2\sigma_{t}(i) + 2\sigma_{t}(j) + \sigma_{t}(i+j)}
\lambda^{i+j} f\big(c^{i}, c^{j}\big) x \\
\hphantom{\big( c^{i} \triangleright \big(c^{j} \triangleright x\big) \big) f\big(c^{i},
c^{j}\big)}{}
 = (-1)^{\sigma_{t}(i+j)} \lambda^{i+j} f\big(c^{i}, c^{j}\big) x
  = f\big(c^{i}, c^{j}\big) \big(c^{i+j} \triangleright x\big),
\end{gather*}
which proves our assertion and completes the proof.
\end{proof}

\begin{Remark}
If $t_{1} \in \mathcal{U} (C_{n}, C_{2})$ is the trivial map, i.e.\
$t_1 (c^i) = 1$, for any $i = 1, \dots, n-1$, then $2 \,|\,
\sigma_{t_{1}} (j)$, for all $j \geq 0$. Then the crossed systems
$(H_{4}, k[C_{n}], \triangleright, f)$ corresponding to $(t_{1},
\lambda)$, $\lambda^{n} = 1$, are precisely the matched pairs of
\cite[Proposition~4.3]{ABM1}.
\end{Remark}

Next we describe by generators and relations all crossed products
$H_4 \# k[C_n]$ associated to the crossed systems from
Theorem~\ref{sistemecros}. We need to introduce the following notation:
for any $i,j = 1, \dots, n-1$ we shall denote by $j \star i$
the following number
\begin{gather*}
j \star i =
\begin{cases} j+i & \mbox{if} \ \  j+i \leq n,\\
j+i-n & \mbox{if} \ \  j+i > n.
\end{cases}
\end{gather*}

\begin{Corollary} \label{anavers}
Let $k$ be a field of characteristic $\neq 2$ and $n$ a positive
integer. Then $H_4 \#^{\triangleright}_f   k[C_n] \cong H_{4n,
\lambda,   t}$, for some $(t, \lambda) \in {\mathcal C} {\mathcal
S} (n, k)$, where we denote by $H_{4n, \lambda,  t}$ the
$4n$-dimensional Hopf algebra having $\{ d_i,   g d_i,   x d_i,
  gx d_i \, | \, i = 1, \dots, n \} $ as a $k$-basis, the unit $1
= d_n$ and the multiplication is subject to the following
relations
\begin{gather*}
g^2  =  1, \qquad x^2 = 0, \qquad   x g = -g x, \qquad d_i g
= g d_i,
\\
d_i   x  =  (-1)^{\sigma_{t}(i)} \lambda^{i}   x   d_i, \qquad
  d_{j} d_{i} = t(c) \cdots t \big(c^{j-1}\big)   t\big(c^i\big) t\big(c^{i+1}\big)
\cdots t \big(c^{i+j-1}\big)   d_{j \star i}
\end{gather*}
for all $i, j = 1, \dots, n-1$. The coalgebra structure and the
antipode on $H_{4n, \lambda,  t}$ are given by
\begin{gather*}
\Delta(g)  =  g \otimes g, \qquad \Delta(d_i) = d_{i} \otimes
d_{i}, \qquad \Delta(x) = x \otimes 1 + g \otimes x, \qquad
\varepsilon(g) =
\varepsilon(d_{i}) = 1, \\
\varepsilon(x)  =  0, \qquad S(g) = g, \qquad S(x) = - g x, \qquad
S(d_{i}) = d_{n-i}   \prod_{p=0}^{i-1}   t\big(c^{p}\big)
t\big(c^{n-1-p}\big)
\end{gather*}
for all $i = 1, \dots, n-1$.
\end{Corollary}

\begin{proof}
The Hopf algebra $H_{4n, \lambda,  t}$ is the crossed product
$H_{4}   \#^{\triangleright}_{f}   k[C_{n}]$ associated to the
pair $(t, \lambda)$ as in Theorem~\ref{sistemecros}. Up to canonical
identif\/ication, the crossed product $H_{4}
\#^{\triangleright}_{f}   k[C_{n}]$ is generated as an algebra by
$g = g \#  1$, $x = x \#  1$ and $d_{i} = 1 \#
c^{i}$ for $i \in \{1, \dots, n\}$. As $H_{4}$ is a Hopf
subalgebra of $H_{4n, \lambda,  t}$, the relations $g^{2} =
1$, $x^{2} = 0$ and $g x = - x g$ also hold in $H_{4n,   \lambda,
  t}$. Now, in the crossed product $H_{4}
\#^{\triangleright}_{f}   k[C_{n}]$ the following relations hold
\begin{gather*}
d_{i} g  =  \big(1 \#  c^{i}\big)(g \#  1) = c^{i}
\triangleright g   \#
  c^{i} = g \#  c^{i} = (g \# 1) \big(1 \# c^i\big) = g d_{i},\\
d_{i} x  =  \big(1 \#  c^{i}\big)(x \#  1) = c^{i}
\triangleright x\#  c^{i} = (-1)^{\sigma_{t}(i)} \lambda^{i}
x \#  c^{i} =
(-1)^{\sigma_{t}(i)} \lambda^{i} x d_{i},\\
d_{j} d_{i}  =  \big(1 \#  c^j\big)\big(1 \#  c^{i}\big) = f\big(c^j,
c^{i}\big)   \# c^{j+i} = t(c) \cdots t \big(c^{j-1}\big)   t\big(c^i\big)
t\big(c^{i+1}\big) \cdots t \big(c^{i+j-1}\big)   d_{j \star i}
\end{gather*}
for all $i,j = 1, \dots, n-1$. The formula for the antipode
follows from~\eqref{antipod}.
\end{proof}

Next we shall give necessary and suf\/f\/icient conditions for two
Hopf algebras $H_{4n, \lambda,  t}$ and $H_{4n,  \lambda',
  t'}$ to be isomorphic. We recall from \cite[Lemma~4.6]{ABM1}
that $\Aut_{\rm Hopf}(H_4) \cong k^*$: explicitly, any
automorphism $u : H_4 \to H_4$ is of the form
\begin{gather}
u (g) = g, \qquad u (x) = \beta   x, \qquad u (gx) = \beta   gx
\label{impleh4}
\end{gather}
for some non-zero scalar $\beta \in k^*$. It what follows, the
automorphism $u$ of $H_4$ implemented by $\beta \in k^*$ as in~\eqref{impleh4} will be denoted by $u = u_{\beta}$.

\begin{Theorem} \label{clasificareH_4C_n}
Let $k$ be a field of characteristic ${\neq} 2, n$ a positive
integer and $(t, \lambda)$, $(t', \lambda') \in {\mathcal C}
{\mathcal S} (n, k)$. Then there is a bijective correspondence
between the set of all Hopf algebra isomorphisms $\psi: H_{4n,
\lambda,   t} \to H_{4n,   \lambda',   t'}$ and the set of all
triples $(\beta,   r,   s) \in k^* \times
\mathcal{U}(C_{n},C_{2}) \times \{ s \in \{1, \dots, n-1 \} \, |
\,   \gcd (s, n) = 1 \}$ satisfying the following
compatibility conditions for any $i,j = 1, \dots, n$
\begin{gather}
\prod_{p=0}^{i-1} t\big(c^{p}\big)   t\big(c^{p+j}\big)  =  r\big(c^{i}\big)   r\big(c^{j}\big)
  r\big(c^{i+j}\big)   \prod_{q=0}^{is - 1}
t'\big(c^{q}\big)   t'\big(c^{q + js}\big), \label{CP6HC} \\
(-1)^{1 + \sigma_{t}(i) + \textnormal{ord} \, r (c^{i})}
\lambda^{i}   \beta  =  (-1)^{\sigma_{t'}(is)}   (\lambda
')^{is} .\label{CP7HC}
\end{gather}

Under the above bijection the isomorphism $\psi = \psi_{(\beta, r,
s)} : H_{4n, \lambda,  t} \to H_{4n,  \lambda',   t'}$
corresponding to $(\beta, r, s)$ is given by
\begin{gather*}
\psi \big(a \# c^{i}\big) = u_{\beta}(a)   r\big(c^{i}\big) \# c^{is}
\end{gather*}
for all $a \in H_{4}$ and $i = 1,\dots, n$.
\end{Theorem}

\begin{proof}
Let $(H_{4}, k[C_n],   \triangleright,   f_t)$ and $(H_{4},
k[C_n],   \triangleright',   f_{t'})$ be the crossed systems
corresponding to $(t, \lambda)$ and respectively $(t',\lambda')$
given in Theorem~\ref{sistemecros}. Then $H_{4n, \lambda,  t} =
H_{4} \#^{\triangleright}_{f_t}   k[C_{n}]$ and $H_{4n,
\lambda',   t'} = H_{4} \#^{\triangleright'}_{f_{t'}}
k[C_{n}]$. By \cite[Theorem~2.1]{ABM2}, the set of all Hopf
algebra morphisms $\psi: H_{4n, \lambda,  t} \to H_{4n,
\lambda',   t'} $ is in bijective correspondence with the set of
all quadruples $(u, p, r, v)$, where $p: H_{4} \to k[C_{n}]$ is a
Hopf algebra map, $u: H_{4} \to H_{4}$, $r: k[C_{n}] \to H_{4}$
and $v: k[C_{n}] \to k[C_{n}]$ are unitary coalgebra maps
satisfying the following compatibility conditions:
\begin{alignat*}{3}
&{\rm (CP1)} \quad && u(a_{(1)}) \otimes p(a_{(2)}) = u(a_{(2)}) \otimes
p(a_{(1)}),&\\
&{\rm (CP2)} \quad && r(h_{(1)}) \otimes v(h_{(2)}) =
r(h_{(2)}) \otimes v(h_{(1)}),&\\
&{\rm (CP3)} \quad && u(ab) = u(a_{(1)})
\big( p (a_{(2)}) \triangleright' u(b_{(1)}) \big) f'
\big(p(a_{(3)}),   p(b_{(2)})\big),&\\
&{\rm (CP4)} \quad && v(h)   v(g)
= p\big(f(h_{(1)},   g_{(1)})\big)   v(h_{(2)}g_{(2)}),&\\
&{\rm (CP5)} \quad && v(h)   p(a) = p(h_{(1)} \triangleright a)
v(h_{(2)}),&\\
&{\rm (CP6)} \quad && r(h_{(1)}) \big(v(h_{(2)})
\triangleright' r(g_{(1)})\big)
  f'   \big(v(h_{(3)}),   v(g_{(2)})\big) & \\
&&& {}= u \big(f(h_{(1)},   g_{(1)})\big) \big(p\big(f(h_{(2)},
g_{(2)}\big) \triangleright' r(h_{(4)} g_{(4)})\big) f'
\big(p\big(f(h_{(3)}, g_{(3)}),  v(h_{(5)}g_{(5)})\big)\big),&\\
&{\rm (CP7)} \quad && r(h_{(1)}) \big(v(h_{(2)}) \triangleright'
u(a_{(1)})\big)
  f'   \big(v(h_{(3)}),   p(a_{(2)})\big) & \\
&&& {}=u(h_{(1)} \triangleright a_{(1)})   \big(p(h_{(2)}
\triangleright a_{(2)}) \triangleright' r(h_{(4)})\big) f'
\big(p(h_{(3)} \triangleright a_{(3)}),   v(h_{(5)})\big)  &
\end{alignat*}
for all $a, b \in H_4$, $g, h \in k[C_n]$. The correspondence
is such that the morphism $\psi = \psi_{(u, p, r, v)}$ associated
to $(u, p, r, v)$ is given by
\begin{gather*}
\psi(a \# h) = u(a_{(1)})   \big( p(a_{(2)}) \triangleright'
r(h_{(1)}) \big)f'   \big(p(a_{(3)}),   v(h_{(2)})\big)
\#'   p(a_{(4)})   v(h_{(3)})
\end{gather*}
for all $a \in H_{4}$ and all $h \in k[C_{n}]$. We will show that,
under this bijection, isomorphisms correspond precisely to
quadruples $(u, p, r, v)$, where $p: H_4 \to k[C_n]$ is the
trivial morphism, $u: H_{4} \to H_{4}$ and $v: k[C_{n}] \to
k[C_{n}]$ are Hopf algebra automorphisms and $r: k[C_{n}] \to
H_{4}$ is a~unitary coalgebra map such that the following two
conditions are satisf\/ied
\begin{gather}
f_{t}\big(c^{i}, c^{j}\big) = r\big(c^{i}\big) r\big(c^{j}\big) r\big(c^{i+j}\big) f_{t'}\big(
v\big(c^{i}\big),v\big(c^{j}\big) \big), \label{CP6HC'}
\\
u\big(c^{i} \triangleright a\big) r\big(c^{i}\big) = r\big(c^{i}\big) \big( v\big(c^{i}\big)
\triangleright' u(a) \big) \label{CP7HC'}
\end{gather}
for all $i = 1,\dots, n$ and $a \in H_{4}$. At the end we will see
that such quadruples are in bijection with triples $(\beta,   r,
  s) \in k^* \times \mathcal{U}(C_{n},C_{2}) \times \{ s \in \{1,
\dots, n-1 \} \, | \,   \gcd (s, n) = 1 \}$ such that~\eqref{CP6HC} and~\eqref{CP7HC} are satisf\/ied.

Suppose f\/irst that $\psi = \psi_{(u, p, r, v)} : H_{4}
\#^{\triangleright}_{f_t} k[C_{n}] \to H_{4}
\#^{\triangleright'}_{f_{t'}} k[C_{n}]$ is an isomorphism
correspon\-ding to $(u, p, r, v)$. Thus $p: H_{4} \to k[C_{n}]$ is a
Hopf algebra map; it follows from \cite[Lemma~4.6]{ABM1} that
$p(x) = 0$ and $p(gx) = 0$. In particular, we have
\begin{gather*}
\psi (x \# 1) = u(x) \# p(1) + u(g) \# p(x) = u(x) \# 1.
\end{gather*}
As $\psi$ is an isomorphism, $u(x)$ must be non-zero. Looking at
(CP1) for $a = x$ and taking into account that $p(x) = 0$ and
$u(x) \neq 0$ we obtain $p(g) = 1$. Thus, $p: H_4 \to k[C_n]$ is
the trivial morphism: $p (z) = \varepsilon (z) 1$, for all $z \in
H_4$. Thus, the compatibility conditions~(CP1) and~(CP5) hold
automatically while (CP2) holds since $k[C_n]$ is cocommutative.
Moreover, (CP3)~is equivalent to the fact that $u : H_{4} \to
H_{4}$ is a Hopf algebra morphism and~(CP4) is equivalent to~$v :
k[C_{n}] \to k[C_{n}]$ being a Hopf algebra morphism. Since $u(x)
\neq 0$, it follows from \cite[Lemma~4.6]{ABM1} that $u: H_4 \to
H_4$ is an automorphism of $H_4$. Using \cite[Corollary~2.2]{ABM2}
we obtain that $v: k[C_n] \to k[C_n]$ is also an automorphism of
$k[C_n]$.

\looseness=-1
It remains to show that~\eqref{CP6HC'} and~\eqref{CP7HC'} hold.
We claim that, in fact, these are exactly (CP6) and (CP7) written
in an equivalent form. Indeed, taking into account that~$p$ is the
trivial map, we have that (CP7) is equivalent to the compatibility
condition~\eqref{CP7HC'}, while (CP6) is equivalent to
\begin{gather*}
r\big(c^{i}\big) \big( v\big(c^{i}\big) \triangleright' r\big(c^{j}\big) \big) f_{t'}
\big( v\big(c^{i}\big), v\big(c^{j}\big) \big) = u \big( f_{t}\big(c^{i}, c^{j}\big) \big)
r\big(c^{i+j}\big)
\end{gather*}
for all $i, j \in \{1, \dots, n\}$. Now $r : k[C_{n}] \to H_{4}$
is a coalgebra map, hence $r(c^{i}) \in \{1, g\}$, for all~$i$.
Since $v(c^{i}) \in C_{n}$ and the elements of~$C_n$ act, via
$\triangleright'$, trivially on $\{1, g\}$, we have $v(c^{i})
\triangleright' r(c^{j}) = r(c^{j})$, for all $i$ and $j$.
Furthermore, $f_{t}(c^{i}, c^{j}) \in \{1,g\}$ for all $i$ and
$j$, and $u(1) = 1$ and $u(g) = g$, hence, $u ( f(c^{i},
c^{j}) ) = f(c^{i}, c^{j})$, for all $i$ and $j$. These
remarks show that (CP7) is equivalent to
\begin{gather*}
r\big(c^{i}\big)   r\big(c^{j}\big)   f_{t'} \big( v\big(c^{i}\big), v\big(c^{j}\big) \big) =
f_{t}\big(c^{i}, c^{j}\big)   r\big(c^{i+j}\big)
\end{gather*}
for all $i,j = 1, \dots, n$, which is equivalent to~\eqref{CP6HC'}.

Conversely, let $(u, p, r, v)$ be a quadruple with $p : H_{4} \to
k[C_n]$ the trivial morphism, $u: H_{4} \to H_{4}$ and $v:
k[C_{n}] \to k[C_{n}]$ Hopf algebra automorphisms and $r: k[C_{n}]
\to H_{4}$ a unitary coalgebra map satisfying~\eqref{CP6HC'} and~\eqref{CP7HC'}. We will prove that the compatibility conditions
(CP1)--(CP7) are satisf\/ied. Indeed, (CP1)--(CP5) are trivially
fulf\/illed. (CP6)~and~(CP7) are equivalent, as we have seen, with~\eqref{CP6HC'} and~\eqref{CP7HC'}. Thus, $(u, p, r, v)$
determines a Hopf algebra morphism $\psi: H_{4}
\#^{\triangleright}_{f_t} k[C_{n}] \to H_{4}
\#^{\triangleright'}_{f_{t'}} k[C_{n}]$, given by
\begin{gather*}
\psi(a \# h) = u(a) r(h_{(1)}) \#' v(h_{(2)})
\end{gather*}
for all $a \in H_{4}$ and $h \in k[C_{n}]$. Since $u$ and $v$ are
isomorphisms it follows from \cite[Corollary~2.2]{ABM2} that
$\psi$ is an isomorphism.

To conclude, we have established a bijective correspondence
between the set of all Hopf algebra isomorphisms $\psi: H_{4n,
\lambda,   t} \to H_{4n,   \lambda',   t'} $ and the set of all
quadruples $(u, p ,r, v)$ consisting of the trivial morphism $p :
H_{4} \to k[C_n]$, two Hopf algebra automorphisms, $u: H_{4} \to
H_{4}$ and $v: k[C_{n}] \to k[C_{n}]$, and a unitary coalgebra map
$r: k[C_{n}] \to H_{4}$ that satisfy~\eqref{CP6HC'} and~\eqref{CP7HC'}. These quadruples are, in turn, in bijection with
the set of all triples $(\beta,   r,   s) \in k^* \times
\mathcal{U}(C_{n},C_{2}) \times \{ s \in \{1, \dots, n-1 \} \, |
\,   \gcd (s, n) = 1 \}$ such that the compatibility conditions~\eqref{CP6HC} and~\eqref{CP7HC} are fulf\/illed. Indeed, for any
$u \in \Aut_{\rm Hopf}(H_{4})$ there exists a unique $\beta \in
k^*$ such that $u = u_{\beta}$ and for any $v \in \Aut_{\rm Hopf}
\big( k[C_{n}] \big)$ there exists a unique $s \in \{1, \dots, n-1
\}$, $\gcd (s, n) = 1$, such that $v(c^{i}) = c^{is}$, for all $i
= 0, \dots, n-1$. Furthermore, any unitary coalgebra map $r:
k[C_{n}] \to H_{4}$ is uniquely determined by a unitary map $C_{n}
\to C_{2} = \{1, g\}$ which we still denote by $r$. Therefore~\eqref{CP6HC} and~\eqref{CP7HC} are nothing but~\eqref{CP6HC'}
and~\eqref{CP7HC'} in terms of $t$, $t'$, $\beta$, $r$ and $s$.
\end{proof}

\begin{Remark} \label{remarca utila}
Although condition~\eqref{CP6HC} is given in terms of $t$, $t'$,
$r$ and $s$, in practice it is rather dif\/f\/icult to work with. In
such a situation, it is more convenient to consider the equivalent
condition~\eqref{CP6HC'} which says that the normalized
2-cocycles $f_{t}$ and $f_{t'} \circ (v \times v)$, where $v$ is
the automorphism of $C_{n}$ associated to~$s$, are cohomologous
and that the coboundary by which they dif\/fer is the one associated
to~$r$.
\end{Remark}

The main result of the paper is the following theorem:

\begin{Theorem} \label{aldoileaseint}
Let $k$ be a field of characteristic $\neq 2$ and $n$ a positive
integer. Then:
\begin{enumerate}\itemsep=0pt
\item[$1.$] There exists a bijection between $\mathcal{H}^{2} (
k[C_{n}], H_4)$ and the quotient set $\mathcal{CS} (n, k)
/{\approx}$, where $\approx$ is the equivalence relation on
$\mathcal{CS} (n, k)$ defined by: $(t, \lambda) \approx (t',
\lambda')$ if and only if there exists $r \in \mathcal{U}
(C_{n},C_{2})$ such that for all $i,j = 1, \dots, n$
\begin{gather}
\prod_{p=0}^{i-1} t\big(c^{p}\big)   t\big(c^{p+j}\big)  =  r\big(c^{i}\big) r\big(c^{j}\big)
r\big(c^{i+j}\big) \prod_{q=0}^{i - 1} t'\big(c^{q}\big)
t'\big(c^{q + j}\big),\label{CP6HC''} \\
(-1)^{1 + \sigma_{t} (i) + \textnormal{ord} \, r(c^{i})}
\lambda^{i}  =  (-1)^{\sigma_{t'}(i)} (\lambda
')^{i}.\label{CP7HC''}
\end{gather}

\item[$2.$] There exists a bijection between
$\textsc{C}\textsc{r}\textsc{p} ( k[C_{n}], H_4)$ and the quotient
set $\mathcal{CS} (n, k) /{\equiv}$, where $\equiv$ is the
equivalence relation on $\mathcal{CS} (n, k)$ defined by:
$(t,\lambda) \equiv (t', \lambda')$ if and only if there exists
$(\beta, r, s) \in k^* \times \mathcal{U} (C_{n}, C_{2}) \times \{
s \in \{1, \dots, n-1 \} \, |  \, \gcd (s, n) = 1 \}$ satisfying
the compatibility conditions~\eqref{CP6HC} and~\eqref{CP7HC}.
\end{enumerate}
\end{Theorem}

\begin{proof}
It follows from Theorems~\ref{sistemecros} and~\ref{clasificareH_4C_n}.
For $(1)$ we use the fact that the isomorphism $\psi =
\psi_{(\beta, r, s)} : H_{n,  t,  \lambda} \to H_{n,  t',
\lambda'}$ associated to the triple $(\beta, r, s) \in k^* \times
\mathcal{U} (C_{n},C_{2}) \times \{ s \in \{1, \dots, n-1 \} \, |
\,  \gcd (s, n) = 1 \}$ stabilizes $H_4$ if and only if $\beta =
1$ and co-stabilizes $k[C_{n}]$ if and only if $s = 1$.
\end{proof}

Theorem~\ref{clasificareH_4C_n} allows us to give a description of the
automorphisms of $H_{4n, \lambda,  t}$. We denote, in what
follows, by $v_{s}$ the automorphism of $C_{n}$ associated to $s
\in \{1, \dots, n-1\}$, $\gcd (s,n) = 1$, given by $v_{s} (c^{i})
= c^{is}$, for all $i$. We also note that $(\mathcal{U} (C_{n},
C_{2}),   \cdot)$ is an Abelian group with point-wise
multiplication.

\begin{Corollary} \label{autom}
The group $\Aut_{\rm Hopf} (H_{4n, \lambda, t})$ of all Hopf
algebra automorphisms is parameterized by the set of all triples
$(\beta, r, s) \in k^* \times \mathcal{U} (C_{n}, C_{2}) \times \{
s \in \{1, \dots, n-1 \} \, |  \, \gcd (s, n) = 1 \}$ satisfying~\eqref{CP6HC} and~\eqref{CP7HC} with $(\lambda', t') =
(\lambda, t)$. The automorphism of $H_{4n, \lambda, t}$
corresponding to $(\beta, r, s)$ is given for any $a \in H_{4}$
and $i = 1,\dots, n$ by
\begin{gather*}
\psi_{(\beta,   r,   s)} : \ H_{4n, \lambda, t} \to H_{4n,
\lambda, t}, \qquad \psi_{(\beta,   r,   s)} \big(a \# c^{i}\big) =
u_{\beta} (a) r \big(c^{i}\big) \#  v_{s} \big(c^{i}\big).
\end{gather*}
In particular, there exists an embedding
\begin{gather*}
\Aut_{\rm Hopf} (H_{4n, \lambda, t}) \hookrightarrow k^{*} \times
\big( \mathcal{U} (C_{n}, C_{2}) \rtimes_{f} \Aut (C_{n}) \big),
\end{gather*}
where $\mathcal{U} (C_{n}, C_{2}) \rtimes_{f} \Aut (C_{n})$ is the
semidirect product associated to $f: \Aut (C_{n}) \to \Aut  (
\mathcal{U} (C_{n},$ $C_{2})  )$,  $f (v) (p) = p \circ v $, for
all $p \in \mathcal{U} (C_{n}, C_{2})$ and $v \in \Aut (C_{n})$.
\end{Corollary}

\begin{proof}
The f\/irst assertion follows from Theorem~\ref{clasificareH_4C_n}. If
$\psi_{(\alpha,  p ,   q)}$ and $\psi_{(\beta,   r,   s)}$ are
two automorphisms of $H_{4n, \lambda, t}$, then, using the fact
that $u_{\alpha}$ is a homomorphism which acts as the identity on
$\{1, g\}$, we have
\begin{gather*}
\psi_{(\alpha,   p,   q)} \circ \psi_{(\beta,   r,   s)} \big(a \#
c^{i}\big)   =  \psi_{(\alpha,   p,  q)} \big( u_{\beta} (a) r
\big(c^{i}\big) \#  v_{s} \big(c^{i}\big) \big)
  =   u_{\alpha} \big( u_{\beta} (a) r \big(c^{i}\big) \big) p \big(c^{is}\big) \#
  v_{q} v_{s} \big(c^{i}\big)\\
 \hphantom{\psi_{(\alpha,   p,   q)} \circ \psi_{(\beta,   r,   s)} \big(a \#
c^{i}\big)}{}
  =   u_{\alpha \beta} (a) r \big(c^{i}\big) (p \circ v_{s}) \big(c^{i}\big) \#
  v_{qs} \big(c^{i}\big)
  =   \psi_{( \alpha \beta,   r \cdot (p \circ v_{s}),   qs)} \big(a
\# c^{i}\big).
\end{gather*}
This shows that $ \Gamma : \Aut_{\rm Hopf} (H_{4n, \lambda, t})
\to k^{*} \times  ( \mathcal{U} (C_{n}, C_{2}) \rtimes_{f} \Aut
(\mathbb{Z}_{n})  )$, def\/ined by $ \Gamma ( \psi_{(\beta,   r,
  s)} ) := \big( \beta, (r, v_{s})^{-1} \big)$ is a one-to-one
homomorphism of groups. Indeed,
\begin{gather*}
\Gamma \big( \psi_{(\alpha,   p,   q)} \circ \psi_{(\beta,   r,
  s)} \big)   =   \Gamma ( \psi_{( \alpha \beta,   r \cdot (p
\circ v_{s}),   q s)} )
  =   \big( \alpha \beta, \big( r \cdot (p \circ v_{s}),
v_{qs} \big)^{-1} \big) \\
\hphantom{\Gamma \big( \psi_{(\alpha,   p,   q)} \circ \psi_{(\beta,   r,
  s)} \big)}{}
  =   \big( \alpha \beta, \big( (r, v_{s}) (p, v_{q})
\big)^{-1} \big)
  =   \big( \alpha,   (p, v_{q})^{-1} \big)   \big( \beta,
(r, v_{s})^{-1} \big) \\
\hphantom{\Gamma \big( \psi_{(\alpha,   p,   q)} \circ \psi_{(\beta,   r,
  s)} \big)}{}
  =  \Gamma (\psi_{(\alpha,   p,   q)} )   \Gamma
(\psi_{(\beta,   r,   s)} ),
\end{gather*}
which proves our claim and concludes the proof.
\end{proof}

\begin{Example}\label{exemplul n=2}
As it can be easily seen from Theorem~\ref{aldoileaseint}, the
description of the classifying objects $\mathcal{H}^{2} (
k[C_{n}], H_4)$ and $\textsc{Crp} (k[C_{n}], H_4)$ depends on the
arithmetics of the positive integer $n$. We describe below the two
classifying objects for $n = 2$~-- this can serve as a model for
other values of $n$. Let $k$ be a f\/ield of characteristic $\neq
2$. Then:
\begin{enumerate}\itemsep=0pt
\item[1.] If $-1 \notin k^{2}$, then $\mathcal{H}^{2} ( k[C_{2}], H_4)
\cong \textsc{Crp} ( k[C_{2}], H_4) = \{   H_4 \otimes k[C_2]
\}$.

\item[2.] If $-1 = \zeta^{2}$, $\zeta \in k$, then $\mathcal{H}^{2} (
k[C_{2}], H_4) \cong \textsc{Crp} ( k[C_{2}], H_4) = \{   H_4
\otimes k[C_2],    H_{8,   \zeta}   \}$, where $H_{8,
\zeta}$ is the $8$-dimensional Hopf algebra having $\{1, g, x, gx,
d, gd, xd, gxd \}$ as a $k$-basis, with the multiplication subject
to the following relations
\begin{gather*}
g^2 = 1, \qquad x^2 = 0, \qquad d^{2} = g, \qquad x g = -g x, \qquad d
g = g d, \qquad d x = \zeta   x d,
\end{gather*}
and the coalgebra structure such that $g$ and $d$ are group-like
elements and $x$ is $(1, g)$-primitive.
\end{enumerate}

Indeed, $\mathcal{U} (C_{2}, C_{2}) = \{ t_{1}, t_{2} \}$, where
$t_{1} (c) = 1$ and $t_{2} (c) = c$. Then $\sigma_{t_{1}} (2) = 2$
and $\sigma_{t_{2}} (2) = 3$. Thus,
\begin{gather*}
 \mathcal{CS} (2, k) =
\begin{cases}
\{ (t_{1}, 1),   (t_{1}, -1) \} & \textnormal{if}  \ \  {-}1 \notin
k^{2},\\
\{ (t_{1}, 1),   (t_{1}, -1),   (t_{2}, \zeta),   (t_{2},
-\zeta)
\} & \textnormal{if}  \ \ {-}1 = \zeta^{2}, \quad \zeta \in k.
\end{cases}
\end{gather*}
It is easy to see that $(t_{1}, 1) \approx (t_{1}, -1)$ and
$(t_{2}, \zeta) \approx (t_{2}, - \zeta)$, the map $r \in
\mathcal{U} (C_{2}, C_{2})$ satisfying~\eqref{CP6HC''} and~\eqref{CP7HC''}, being, in both cases, $r = t_{2}$. Since the
normalized $2$-cocycles, $f_{t_1}$ and $f_{t_{2}}$, associated to~$t_{1}$ and~$t_{2}$ are not cohomologous, and the~$s$ from~\eqref{CP6HC} can take only the value~1, we obtain using Remark~\ref{remarca utila}, that $(t_{1}, 1) \not \equiv (t_{2},
\zeta)$. Thus, $(t_{1}, 1) \not \approx (t_{2}, \zeta)$, also.
Finally, the description of~$H_{8,   \zeta}$ follows from Corollary~\ref{anavers}.
\end{Example}

\subsection*{Acknowledgements}

 The authors would like to thank the referees for their
comments and suggestions that substantially improved the f\/irst
version of this paper.
A.L.~Agore is research fellow `Aspirant' of FWO-Vlaanderen. This
work was supported by a grant of the Romanian National Authority
for Scientif\/ic Research, CNCS-UEFISCDI, grant no.~88/05.10.2011.

\pdfbookmark[1]{References}{ref}
\LastPageEnding

\end{document}